\documentclass[letterpaper,12pt]{article}
\usepackage{tabularx} 
\usepackage{amsmath}  
\usepackage{enumitem}
\usepackage{longtable}
\usepackage{array}
\usepackage{authblk}
\usepackage{graphicx} 
\usepackage[margin=1in,letterpaper]{geometry} 
\usepackage{cite} 
\usepackage[final]{hyperref} 
\hypersetup{
	colorlinks=true,       
	linkcolor=blue,        
	citecolor=blue,        
	filecolor=magenta,     
	urlcolor=blue         
}
\usepackage{blindtext}

\begin{document}

\title{Predicting Winning Lottery Numbers}
\author{Thando Nkomozake}
\affil[]{Department of Mathematics and Physics, Cape Peninsula University of Technology,
Bellville, Cape Town, South Africa}

\maketitle

\begin{abstract}
\noindent We use mathematical statistics theory to derive the Compound-Dirichlet-Multinomial (CDM) prediction model. We then use this model to predict winning numbers for the 6-number, 5-number, pick-4 and pick-3 lottery games. We also develop a strategy which we call ``the 3-strategy", for generating profit by predicting winning numbers for the pick-3 lottery game.
\end{abstract}

\section{Introduction}
\noindent The purpose of this work is to obtain a good model for predicting winning lottery numbers. In chapter 2 we discuss the prior and posterior distributions in some detail. We then provide a practical example based on the Bernoulli distribution. In chapter 3 we focus on the vector case of our theoretical discussion. In section 3.1 we discuss the definition of the Multinomial distribution. In section 3.2 we discuss the definition of the Dirichlet distribution and use this and our results in the previous sections to define the Compound-Dirichlet-Multinomial distribution.\\

\noindent In chapter 4 we focus on the matrix case of our theoretical discussion. In section 4.1 we derive the probability distribution function (pdf) of the Compound-Dirichlet-Multinomial distribution for the matrix case. In section 4.2 we use the results of section 4.1 to write the (CDM) prediction model. In chapter 5 we discuss the three methods that we use to estimate the parameter $ \overline{\alpha} $ of the (CDM) prediction model. In section 5.1 we discuss maximum likelihood estimation, in section 5.2 we discuss the method of moments and in section 5.3 we discuss the main diagonal method. In chapter 6 we discuss our results, where in section 6.1 we discuss our results for the 6-number lottery game and in section 6.2 we discuss our results for the 5-number lottery game. In section 6.3 we discuss our results for the pick-3 lottery game. We introduce our strategy that we call ``the 3-strategy" where we use the (CDM) prediction model to obtain profit from the Vermont pick-3 lottery game. We don't discuss the pick-4 lottery game because our results for it are similar to our results for the pick-3 lottery game.

\section{The Prior and Posterior Distributions}

\subsection{Prior Distribution}

\noindent In the frequentist paradigm the uknown parameter is a fixed number. In the Bayesian paradigm the uknown parameter is a random variable. 
For two random variables $X$ and $Y$, their joint pdf is denoted with $f_{X,Y}(x,y)$. Then the marginal pdf of $X$ is given by, \begin{equation}
    f_{X}(x) = \int f_{X,Y}(x,y) \ dy.
\end{equation}

\noindent The conditional pdf of $Y$ given that $X = x$ is given by, \begin{equation}
    f_{Y|X}(y|x) = \frac{f_{X,Y}(x,y)}{f_{X}(x)}.   
\end{equation}

\noindent In Bayesian analysis the uknown parameter is represented by a random variable $\Theta$. The corresponding pdf is denoted by $f_{\Theta}(\theta)$ and is called the prior distribution. It is given such a name because it gives us information about the uknown parameter prior to the observation of data. 

\subsection{Posterior Distribution}

We now observe the data and denote the observed data with the random vector $\overline{X}$. The conditional pdf of $\overline{X}$ given $\Theta = \theta$ is given by, \begin{equation}
    f_{\overline{X}|\Theta}(\overline{x}|\theta) = \frac{f_{\overline{X},\Theta}(\overline{x},\theta)}{ f_{\Theta}(\theta) }.
\end{equation}

\noindent You can rewrite (3) as, \begin{equation}
    f_{\overline{X},\Theta}(\overline{x},\theta) = f_{\overline{X}|\Theta}(\overline{x}|\theta) \ f_{\Theta}(\theta).
\end{equation} 

\noindent The marginal pdf of $\overline{X}$ in this case is given by, \begin{equation}
    f_{\overline{X}}(\overline{x}) = \int f_{\overline{X}, \Theta} (\overline{x}, \theta) \ d\theta.
\end{equation}

\noindent Then the conditional pdf of $\Theta$ given that $\overline{X} = \overline{x}$ is given by, \begin{equation}
    f_{\Theta|\overline{X}}(\theta|\overline{x}) = \frac{f_{\overline{X}, \Theta} (\overline{x}, \theta)    }{ \int f_{\overline{X}, \Theta} (\overline{x}, \theta) \ d\theta     }.
\end{equation}

\noindent This pdf (6) is called the posterior distribution because it gives us information about the uknown parameter after we observe the data. We can rewrite (6) as follows, \begin{equation}
   f_{\Theta|\overline{X}}(\theta|\overline{x}) = \frac{f_{\overline{X}|\Theta}(\overline{x}|\theta) \ f_{\Theta}(\theta)     }{ \int f_{\overline{X}|\Theta}(\overline{x}|\theta') \ f_{\Theta}(\theta') \ d \theta'  }.
\end{equation} 

\noindent Then we summarize (7) as follows, \begin{equation}
    f_{\Theta|\overline{X}}(\theta|\overline{x}) \propto f_{\overline{X}|\Theta}(\overline{x}|\theta) \ f_{\Theta}(\theta).     
\end{equation}

\subsection{Example}

Let $X \equiv $ getting a heads when flipping a loaded coin. Therefore $X \sim $ Bernoulli($p$), where $p$ is the parameter of the Bernoulli distribution and it denotes the probability of getting a heads. We have that $P \sim $ Uniform$(0, 1)$ and it has the following pdf, \begin{equation}
    f_{P}(p) = 1, \ \forall p \in (0, 1).
\end{equation}

\noindent The pdf of $X$ given that $P = p$ is, \begin{equation}
    f_{X| P}(x| p) = p^{x}  \left(  1 - p  \right)^{1-x} ,
\end{equation}

\noindent where $x \in \{0, 1 \}$. If we have a random vector $\overline{X} = \left( X_{1}, X_{2},...,X_{n} \right)$ that follows the Bernoulli($p$) distribution then we have that, \begin{equation}
\begin{split}
    f_{\overline{X}|P}(\overline{x}|p) & = \prod_{i = 1}^{n} \ p^{x_{i}}  \left( 1 - p\right)^{1 - x_{i}}  \\
    & = p^{\sum_{i = 1}^{n} x_{i}} \left(  1 - p\right)^{n - \sum_{i = 1}^{n} x_{i}} 
    \end{split} .
\end{equation} 
\noindent If we let $\sum_{i = 1}^{n} x_{i} = s$, then we can write (11) as, \begin{equation}
    f_{\overline{X}|P}(\overline{x}|p) = p^{s}\ \left(1 - p  \right)^{n - s} .
\end{equation}
\noindent Using (4) we have that, \begin{equation}
\begin{split} 
f_{\overline{X},P}(\overline{x}, p) & = f_{\overline{X}|P}(\overline{x}|p) \ f_{P}(p)\\
& = p^{s} \left( 1 - p \right)^{n - s}.
\end{split}
\end{equation} 

\noindent Now we can calculate the marginal pdf of $\overline{X} $, \begin{equation}
    f_{X}(\overline{x}) = \int_{0}^{1} \ p^{s} \left( 1 - p \right)^{n - s} \ dp.
\end{equation}

\noindent Therefore the posterior distribution of $P$ is given by, \begin{equation}
    f_{P|\overline{X}} (p|\overline{x}) = \frac{p^{s} \left( 1 - p \right)^{n - s}           }{   \int_{0}^{1} \ p^{s} \left( 1 - p \right)^{n - s} \ dp              }.
\end{equation}

\newpage

\section{The Vector Case}

\subsection{The Multinomial Distribution}

\noindent Examples of categorical variables include the following:
\begin{enumerate}[label=(\alph*)]
        \item Gender has two categories $(k = 2)$ and these are male and female.
        \item The outcome of flipping a six-sided die has six categories $(k = 6)$ and these are 1, 2, 3, 4, 5 and 6.
        
    \end{enumerate}

\noindent The multinomial distribution is commonly used to characterize categorical random variables. Suppose that $Z$ is a categorical random variable with $k$ categories. Let $P(Z = j) = p_{j} \ \forall j = 1, 2,...,k. $ The parameter $\overline{p} = (p_{1}, p_{2},..., p_{k})$ describes the entire distribution of $Z$, where $\sum_{j=1}^{k} p_{j} = 1 $. We now generate a random sample $Z_{1}, Z_{2},..., Z_{n}$. Let, \begin{equation}
    X_{j} = \sum_{i=1}^{n} I \ \{Z = j\}.
\end{equation} 

\noindent So then we have that the random vector $\overline{X} = \left( X_{1}, X_{2},..., X_{k}  \right)      $ follows the multinomial distribution with parameters $(n ; p_{1}, p_{2},..., p_{k})$. We write $\overline{X} \sim MN_{k} \left(n; p_{1}, p_{2},..., p_{k} \right) $. Therefore the pdf of $\overline{X}$ given that $\overline{P} = \overline{p}$ is, \begin{equation}
    f_{\overline{X}|\overline{P}}(\overline{x}|\overline{p}) = \frac{n!      
 }{ x_{1}! \ x_{2}! \ ... \ x_{k}!              } \ \prod_{j = 1}^{k} \ p_{j}^{x_{j}}.
\end{equation}

\subsection{The Dirichlet Distribution}

The Dirichlet distribution is the prior distribution to the multinomial parameters $\overline{p} = \left( p_{1}, p_{2},..., p_{k}  \right)  $. Therefore $ \overline{p} \sim Dir(\overline{\alpha})  $, which implies that its pdf is given by, \begin{equation}
    f_{\overline{P}}((\overline{p})) = \frac{ 
  \Gamma(\sum_{i = 1}^{k} \alpha_{i})  }{ \prod_{i = 1}^{k} \Gamma(\alpha_{i})         } \ \prod_{i = 1}^{k} \ p_{i}^{\alpha_{i} - 1} ,
\end{equation} where $ \alpha_{i} > 0 \ \ \forall  i = 1, 2,..., k. $ \\

\noindent The posterior distribution of $\overline{P}$ is given by, \begin{equation}
\begin{split} 
f_{\overline{P}| \overline{X} }(\overline{p}|\overline{x} ) & \propto f_{\overline{X}| \overline{P}}(\overline{x}|\overline{p}) \ f_{\overline{P}}(\overline{p})              \\
& = \frac{ n! \ \Gamma( \sum_{i = 1}^{k} \alpha_{i}    )         }{ \prod_{i = 1}^{k} x_{i}! \  \prod_{i =1}^{k} \Gamma 
(\alpha_{i})} \ \prod_{i = 1}^{k} \ p_{i}^{(x_{i} + \alpha_{i}) - 1} \\
& \propto  \prod_{i = 1}^{k} \ p_{i}^{(x_{i} + \alpha_{i})-1} \\
& \sim Dir(x_{1} + \alpha_{1}, x_{2} + \alpha_{2},..., x_{k} + \alpha_{k} ).
\end{split}
\end{equation} \\
We also have that the joint pdf of $\overline{X}$ and $\overline{P}$ is given by,  \begin{equation}
\begin{split}
    f_{\overline{X}, \overline{P}}(\overline{x}, \overline{p}) & = f_{\overline{X}|\overline{P}}(\overline{x}| \overline{p}) \ f_{\overline{P}}(\overline{p})  \\
    & = \frac{ n! \ \Gamma(\sum_{i = 1}^{k} \alpha_{i})   }{ \prod_{i = 1}^{k} x_{i}! \ \prod_{i = 1}^{k} \Gamma(\alpha_{i})           } \ \prod_{i = 1}^{k} \ p_{i}^{(x_{i} + \alpha_{i}) - 1   } .
    \end{split} 
\end{equation}  \\
The marginal pdf of $\overline{X}$ is the integral of the joint pdf (20), so we have, \begin{equation}
\begin{split}
   f_{\overline{X}}(\overline{x})  & = \int \frac{ n! \ \Gamma(\sum_{i = 1}^{k} \alpha_{i}) }{ \prod_{i = 1}^{k} x_{i}! \  \prod_{i = 1}^{k}  \Gamma(\alpha_{i})       } \ \prod_{i = 1}^{k} \ p_{i}^{( x_{i} + \alpha_{i}) - 1 }  \ d \overline{p} \\                 
    & =  \frac{n! \ \Gamma(\sum_{i = 1}^{k} \alpha_{i}) \ \prod_{i = 1}^{k} \Gamma(x_{i} + \alpha_{i})    }{ \prod_{i = 1}^{k} x_{i}! \ \prod_{i = 1}^{k} \Gamma(\alpha_{i}) \ \Gamma( \sum_{i = 1}^{k} (x_{i} + \alpha_{i})    )          }  .\end{split} 
\end{equation}    \\
This marginal pdf of $\overline{X}$ (21) is called the Compound-Dirichlet-Multinomial (CDM) distribution. Its corresponding expectation value is given by, \begin{equation}
    E(x_{i}) = n \ \frac{ \alpha_{i} }{  \sum_{i = 1}^{k} \alpha_{i}    } .
\end{equation}

\section{The Matrix Case}

\subsection{The (CDM) Distribution For a Matrix}

\noindent Let $ X = \left[ \overline{X}_{i}    \right] \ \forall i = 1, 2,..., n.  $ Where $ \overline{X}_{i} = \left( X_{i1} , X_{i2},..., X_{iK}  \right) $. So $X$ is a matrix with $n$ rows and $K$ columns. Let $ \overline{Y} = \left( y_{1}, y_{2},...y_{K} \right) $ be the multinomial parameter. Let $ M = \sum_{j = 1}^{K} X_{ij} $ and $ n_{j} = \sum_{i = 1}^{n} X_{ij} $. We then have that $ X \sim MN_{K} \left( M ; \overline{Y}  \right) $ and $ \overline{Y} \sim Dir(\overline{\alpha}) $, which implies that the corresponding marginal pdf is given by, \begin{equation}
      \begin{split}
    f_{X}(x) &= \int  \frac{M! \ \Gamma (\sum_{j = 1}^{K} \alpha_{j}  )      }{  \prod_{j = 1}^{K} n_{j} ! \ \prod_{j = 1}^{K} \Gamma(\alpha_{j}) } \ \prod_{j = 1}^{K} \ y_{j}^{ ( n_{j} + \alpha_{j}     ) - 1 }  \ d \overline{y}   \\
    &= \frac{ M! \ \Gamma( \sum_{j = 1}^{K} \alpha_{j} ) \ \prod_{j = 1}^{K} \Gamma(n_{j} + \alpha_{j})   }{  \prod_{j = 1}^{K} n_{j}! \ \prod_{j = 1}^{K} \Gamma(\alpha_{j}) \ \Gamma( \sum_{j = 1}^{K} (n_{j} + \alpha_{j}  ) )   } .
    \end{split}
\end{equation}

\noindent Therefore this pdf (23) is the (CDM) distribution for the matrix $X$. Where, $ \sum_{j = 1}^{K} \alpha_{j} = \alpha_{0} $ and $  \sum_{j = 1}^{K} y_{j} = 1 $. The corresponding expectation value for this distribution is, \begin{equation}
    E(n_{j}) = M \ \frac{ \alpha_{j} }{ \alpha_{0}  } .
\end{equation}    

\subsection{The (CDM) Prediction Model}

\noindent We start by making the following substitution in (23), $ n_{j} \rightarrow z_{j}  $ and $ \alpha_{j} \rightarrow \alpha_{j} + n_{j} $. This gives, \begin{equation}
    f_{Z}(z) = \frac{ M! \ \Gamma( \sum_{j = 1}^{K}( \alpha_{j} + n_{j}   )) \ \prod_{j = 1}^{K} \Gamma(z_{j} + n_{j} + \alpha_{j})          }{ \prod_{j = 1}^{K} z_{j}! \ \prod_{j = 1}^{K} \Gamma(\alpha_{j} + n_{j}) \  \Gamma(\sum_{j = 1}^{K}(z_{j} + n_{j} + \alpha_{j}))              },
\end{equation}

\noindent where $Z$ represents the predicted matrix. The corresponding expectation value is given by, \begin{equation}
    E(z_{j}) = M \ \frac{ (\alpha_{j} + n_{j} )  }{ \sum_{j = 1}^{K} (\alpha_{j} + n_{j}) }.
\end{equation}
Therefore if we want to predict the next row of the matrix $X$ using its previous $n$ rows we use, \begin{equation}
    Pred(x_{n+1}) = M \ \frac{ (\alpha_{j} + n_{j} ) }{ \sum_{j = 1}^{K} (\alpha_{j} + n_{j}) },
\end{equation}
and this is the (CDM) prediction model.

\section{Parameter Estimation}

\subsection{Maximum Likelihood Estimation}

\noindent Suppose that we have the following matrix, \\

\[
\qquad
X = \begin{bmatrix} 
    x_{11} & x_{12} &\dots  & x_{1K}\\
    x_{21} & x_{22} &\dots  & x_{2K}\\
    \vdots & \vdots & \ddots &  \vdots \\
    x_{n1} & x_{n2} &\dots  & x_{nK}
    \end{bmatrix},
\]
where $M = \sum_{j = 1}^{K} x_{ij} $ and $ n_{j} = \sum_{i = 1}^{n} x_{ij} $. Suppose that we want to predict the $ (n + 1)$-th row from the matrix $X$ using the (CDM) prediction model (27). We therefore first have to estimate the parameter $ \overline{\alpha} = (\alpha_{1}, \alpha_{2},..., \alpha_{K}) $. In order to obtain the maximum likelihood estimate for $ \overline{\alpha}  $ we use the following formula [1], \begin{equation}
    \hat{\alpha}_{j} = \alpha_{0} \ f_{j}, \ \forall j = 1, 2,..., K,
\end{equation} where $ f_{j} =  (\sum_{i = 1}^{n} x_{ij})/n   $, and $\alpha_{0}$ is given by the following formula, \begin{equation}
    \alpha_{0} = \frac{ n \ (K - 1) \ \gamma  }{ n \ \sum_{j = 1}^{K} f_{j} \ln({f_{j}}) - \sum_{j = 1}^{K} f_{j} \ \sum_{i = 1}^{n} \ln({x_{ij}})       },
\end{equation} where $\gamma$ is the Euler-Mascheroni constant and is given by $ \gamma = 0.57721566490 $. 

\subsection{Method of Moments}

\noindent The method of moments estimate for the parameter $ 
\overline{\alpha} $ is given by [2], \begin{equation}
    \hat{\alpha}_{j} = \frac{ n_{j}    }{  n    }, \  \forall j = 1, 2,..., K.
\end{equation}

\subsection{Main Diagonal Method}

\noindent The author found that if $X$ is a square matrix then the main diagonal of $X$ is a good estimate for the parameter $\overline{\alpha}$. That gives, \begin{equation}
    \hat{\alpha}_{j} = x_{jj}, \  \forall j = 1, 2,..., K.
\end{equation}

\section{Results}

\subsection{The Six-Number Lottery Game}

\noindent The author collected historical data from the South African national lottery for the six-number lottery game. To play the South African six-number lottery game you pick six numbers from 1 to 52. The historical data shows the past winning numbers for 1971 draws $\equiv$ 21 years. The author then wrote code that compares the numbers predicted by the (CDM) model (27) for each historical draw with the actual winning numbers for that draw. When the model predicts correctly then the code prints out the corresponding winning numbers and their position in the historical data. \\

\noindent We then used the code to compare for 2 winning numbers, 3 winning numbers, 4 winning numbers, 5 winning numbers and 6 winning numbers in the six-number lottery game. We summarized the results on Table 1 below. \\

\begin{table}[ht]
\begin{center}

\label{tbl:bins} 
\begin{tabular}{|c|c|} 
\hline
\multicolumn{1}{|c|}{\textbf{WINNING NUMBERS}} & \multicolumn{1}{c|}{\textbf{TIME DIFFERENCE}} \\
\hline
TWO &   12 draws $\equiv$ \ 1 month \\
THREE &   105 draws $\equiv$ 13 \ months  \\
FOUR &   529 draws $\equiv$ 5.5 \ years \\
FIVE &   9984 draws $\equiv$ 104 \ years \\
SIX &    1000 0032   draws $\equiv$ 104 167 \ years \\
\hline
\end{tabular}
\caption{ (Using the case of two winning numbers to explain. This therefore means that if the code predicts two winning numbers correctly it will take approximately 12 draws $\equiv$ 1 month for it to correctly predict two winning numbers again.) } 
\end{center}
\end{table}

\noindent In using our code on the historical data that covers 21 years, we didn't get any outputs for five winning numbers and six winning numbers. Therefore the time difference between five winning numbers and six winning numbers reflected on Table 1 is approximated using patterns observed in the outputs for two winning numbers, three winning numbers and four winning numbers.

\newpage

\begin{longtable}{|c|c|c|}
\caption{ (These are some outputs from our code. The combinations labelled [MD] are outputs the code prints having used the main diagonal method. The combinations labelled [MM] are outputs the code prints having used the method of moments. The combinations labelled [AC] are the actual winning jackpot numbers for that draw.)  }\\
\hline
(09/07/2022) & (16/07/2022) & (06/08/2022) \\
\hline
\endfirsthead
\caption{\textit{(Continued)} Caption for the multi-page table}\\
\hline
First Column & Second Column & Third Column \\
\hline
\endhead
\hline
\multicolumn{3}{r}{\textit{Continued on next page}} \\
\endfoot
\hline
\endlastfoot
11 \ 19 \ 27 \ 37 \ 39 \ 45 \ [MD] & 8 \ 19 \ 23 \ 31 \ 37 \ 46 \ [MD] & 7 \ 12 \ 27 \ 35 \ 43 \ 47 \ [MD] \\ 11 \ 19 \ 28 \ 36 \ 39 \ 45 \ [MM] & 9 \ 17 \ 23 \ 31 \ 37 \ 46 \ [MM] & 7 \ 13 \ 24 \ 34 \ 45 \ 47 \ [MM] \\ 6 \ 24 \ 29 \ 35 \ 41 \ 44 \ [AC] & 7 \ 13 \ 22 \ 31 \ 45 \ 46 \ [AC] & 7 \ 18 \ 30 \ 36 \ 44 \ 47 \ [AC] \\

\end{longtable}

\subsection{The Five-Number Lottery Game}

\noindent We collected historical data from the Oklahoma Cash 5 lottery which is based in America. The historical data shows the past winning numbers for 4642 draws $\equiv$ 13 years. We then used the same code we used for the six-number lottery game to compare for two winning numbers, three winning numbers, four winning numbers and five winning numbers in the five-number lottery game. We summarized the results on Table 3 below.

\begin{table}[ht]
\begin{center}

\label{tbl:bins} 
\begin{tabular}{|c|c|} 
\hline
\multicolumn{1}{|c|}{\textbf{WINNING NUMBERS}} & \multicolumn{1}{c|}{\textbf{TIME DIFFERENCE}} \\
\hline
TWO &   12 draws $\equiv$  2 \ weeks \\
THREE &   109 draws $\equiv$ 4 \ months  \\
FOUR &   1182 draws $\equiv$ 3 \ years \\
FIVE &   7665 draws $\equiv$ 21 \ years \\

\hline
\end{tabular}
\caption{(Using our code on the historical data from the Oklahoma Cash 5 lottery, that covers 13 years, we didn't get any outputs for five winning numbers. Therefore the time difference between five winning numbers reflected on this table is approximated using patterns observed in the outputs for two winning numbers, three winning numbers and four winning numbers.)}                                                
\end{center}
\end{table}

\subsection{The Pick 3-Lottery Game}

\noindent The author collected historical data from the Vermont Pick 3 lottery, which is based in America. This data contains 5844 draws where 60 draws $\equiv$ 1 month. In order to play the Vermont Pick 3-lottery game you pick three digits between 0 and 9. You win the jackpot if the three digits you picked match the drawn three digits in the exact order. The jackpot prize money is $ \$ 500 $. If there's multiple jackpot winners of a particular draw they each get $ \$ 500  $. There are two draws per day every day of the year. For each draw one person can play up to 21 consecutive three number combinations. The price to play each combination is $ \$ 1  $. We applied our code to the historical data and got the results shown on Table 4 below. Using our results and considering the rules of the game we developed a strategy which we call ``The 3-Strategy" and we discuss it below. With this strategy we play the Vermont Pick 3-lottery game and generate a profit. 

\begin{table}[h]
\centering
\scalebox{1.1}{\begin{tabular}{|c|c|}
    \hline
    \textbf{Draw Number} & \textbf{Difference} \\
    \hline
    0 &  \\
    \hline
     44 & 44 \\
    \hline
     659 & 615 \\
    \hline
     1357 & 698 \\
    \hline
     1369 & 12 \\
    \hline
     1915 & 546 \\
    \hline
     2039 & 124 \\
    \hline
     3449 & 1410 \\
    \hline
     3685 & 236 \\
    \hline
     4285 & 600 \\
    \hline
\end{tabular}}
\caption{(The average difference between two consecutive draws is 476 draws $\equiv$ 8 months. The left column of this table shows the draw numbers where our code correctly predicted the jackpot numbers. The right column shows the difference between two consecutive draw numbers.) }
\label{table-example}
\end{table}

\subsection*{The 3-Strategy: }

\noindent \underline{\textbf{Stream 1}} \\

\noindent The average number of draws between successfully predicted wins is 476 draws $\equiv$ 240 days. We divide 240 days into four quarters and there's 60 days per quarter. We start by using 1 player and we start playing in \underline{the first 60 days:} \begin{equation*}
    \$ 2 \times 60 = \$ 120.
\end{equation*}

\noindent If we win we get $ \$ 500  $ and therefore a profit of $ \$ 380 $. If we don't win in the first 60 days we play in the second 60 days and this time we use 2 seperate players in order to make back the $ \$ 120  $ we lost and also make some profit.\\

\noindent  \underline{The second 60 days:(2 players)}

\begin{equation*}
    \$ 4 \times 60 = \$ 240.
\end{equation*}

\noindent If we win we get $ \$ 1000 $ and therefore a profit of $ \$ 760 $. So $ \$ 760 - \$ 120 = \$ 640 $. We therefore make back our $ \$ 120  $ and make a total profit of $ \$ 640 $. If we don't win in the second 60 days we play in the third 60 days. Also this time we use 5 seperate players in order to make back the $ \$ 120 + \$ 240 = \$ 360 $ we lost and also make some profit.\\

\newpage

\noindent  \underline{The third 60 days:(5 players)}

\begin{equation*}
    \$ 10 \times 60 = \$ 600.
\end{equation*}

\noindent If we win we get $ \$ 2500 $ and therefore a profit of $ \$ 1900 $. So $ \$ 1900 - \$ 360 = \$ 1540 $. We therefore make back our $ \$ 360 $ that we had lost and make a total profit of $ \$ 1540  $. If we don't win in the third 60 days then there's a good probability that we'll win in the fourth 60 days because it completes the 240 days. So we play in the fourth 60 days and this time we use 12 seperate players in order to make back the $  \$ 120 + \$ 240 + \$ 600 = \$ 960  $ we lost and also make some profit.\\

\noindent  \underline{The fourth 60 days:(12 players)}

\begin{equation*}
    \$ 24 \times 60 = \$ 1440.
\end{equation*}

\noindent If we win we get $ \$ 6000 $ and therefore a profit of $ \$ 4560 $. So $ \$ 4560 - \$ 960 = \$ 3600$. We therefore make back our $ \$ 960 $ that we had lost and make a profit of $ \$ 3600 $. That completes the 240 days. In total we spent $ \$ 2400 $ on tickets and won a total of $ \$6000 $ and made a profit of $ \$ 3600 $.\\

\noindent Remember the 240 days is just the average time between successfully predicted wins by our model. Therefore practically this time can be less or more. With respect to the historical data we collected, the most this time has been is 1410 draws $ \equiv $ 730 days. So if it happens that we don't win even in the fourth 60 days, we continue with the procedure, adding more players accordingly until we win. Therefore in this stream 1 of the 3-strategy it's safe to have a budget of $ \$ 100 \ 000   $ in order to play. We can run the 3-strategy across multiple streams simultaneously for maximum profit. If you don't have $ \$ 100 \ 000 $ at your disposal and you want to apply the 3-strategy, you can then observe the pattern displayed by the data and only play at favourable times. The historical data shows that there's ``Short-stretches" and ``Long-stretches". Where a, \begin{equation*}
    \rm{Long-stretch} \geq 500 \ \rm{draws},
\end{equation*}
and a, \begin{equation*}
    \rm{Short-stretch} < 500 \ \rm{draws}.
\end{equation*}

\noindent The data shows that $ 60 \%  $ of the time, the time between successfully predicted wins (by our model) alternates between a Short-stretch and a Long-stretch. Therefore if you don't have $ \$ 100 \ 000  $ in your budget then it's better for you to play during a Short-stretch.

\affil[1]{Department of Mathematics and Physics, Cape Peninsula University of Technology,
Bellville, Cape Town, South Africa}

\end{document}